\documentclass[11pt]{article}
\catcode`\@=11 \catcode`\@=12

\usepackage{amsmath}
\usepackage{amssymb}
\usepackage{amsthm}
\usepackage{oldgerm}

\usepackage[pdftex,colorlinks,urlcolor=blue,pdfstartview=FitH]{hyperref}

\newcommand{\Rm}{\mathbb{R}}

\newcommand{\mP}{\ensuremath{\mathcal{P}}}

\newcommand{\mQ}{\ensuremath{\mathcal{Q}}}

\newcommand{\mD}{\ensuremath{\mathcal{D}}}

\newcommand{\mA}{\ensuremath{\mathcal{A}}}

\newcommand{\Tm}{\ensuremath{\mathbb{T}}}
\newcommand{\e}{\epsilon}
\newcommand{\2}{\frac{1}{2}}

\newtheorem{lem}{Lemma}
\newtheorem{thm}{Theorem}

\newtheorem{hyp}{Hypothesis}

\newtheorem{prop}[lem]{Proposition}
\newtheorem{defn}[lem]{Definition}

\def\proof {\noindent{\sc{Proof. }}}
\def\qed {\mbox{}\hfill {\small \fbox{}} \\}
\def\lto{\longrightarrow}
\def\lmto{\longmapsto}

\def\leq{\leqslant}
\def\geq{\geqslant}

\title{Large normally hyperbolic cylinders in  \textit{a priori} 
stable Hamiltonian systems}

\author{Patrick Bernard\footnote{membre de l'IUF} }

\date{April 2010}

\begin{document}

\maketitle

\begin{small}
\begin{center}
------\\
\noindent
Patrick Bernard,\\
CEREMADE, UMR CNRS 7534\\
Pl. du Mar\'echal de Lattre de Tassigny\\
75775 Paris Cedex 16,
France\\
\texttt{patrick.bernard@ceremade.dauphine.fr}\\
------
\end{center}
\end{small}
\vspace{1cm}
\textbf{Abstract:} We prove the existence of normally hyperbolic cylinders in \textit{a priori}
stable Hamiltonian systems the size of which   is bounded from below independently 
of the size of the perturbation. This result should have applications to the study
of Arnold's diffusion.
\begin{center}
------
\end{center}
\begin{center}
------
\end{center}

A major problem in dynamical systems consists
in studying the Hamiltonian systems on $\Tm^n\times \Rm^n$ of the form

\begin{equation}\tag{$H$}\label{H}
H(q,p)=h(p)-\epsilon^2 G(t,q,p), \quad (t,q,p)\in \Tm\times \Tm^n\times \Rm^n.
\end{equation}
Here $\epsilon$ should be considered as a small perturbation parameter,
we put a square because the sign of the perturbation will play a role
in our discussion.
In  the unperturbed system ($\epsilon =0$)
the momentum variable $p$ is constant.

We want to study the dynamics of the perturbed system in the neighborhood 
of a torus $\{p=p_0\}$, corresponding to a resonant frequency.
There is no loss of generality in assuming that 
the frequency
is of the form 
$$
\partial h(p_0)=(\omega,0)\in \Rm^m \times \Rm^r.
$$
If the restricted frequency $\omega$ is non-resonant in
$\Rm^m$, 
 then it is expected that the averaged system

\begin{equation}\tag{$H_a$}\label{Ha}
H_a(q,p)=H_a(q_1,q_2,p_1,p_2)=h(p)-\epsilon ^2 V(q_2)
 \end{equation}
should locally approximate the dynamics of (\ref{H}) near $p=p_0=(p_1^0,p_2^0)$,
where $q=(q_1,q_2)\in \Tm^m\times \Tm^r$ and $p=(p_1,p_2)\in \Rm^m\times \Rm^r$,
and
where 
$$V(q_2)=\int G(t,q_1,q_2,p_0) dt dq_1.
$$

We  make the following  hypothesis on the averaged system:
\begin{hyp}\label{nondeg}
The function $h$ is convex with positive definite Hessian and
the averaged potential $V$ has a non-degenerate local maximum at $q_2=0$.
\end{hyp}
Under Hypothesis \ref{nondeg},  the  averaged system has an invariant manifold
of equations
$$(\partial_{p_2}h=0, q_2=0)\in \Tm^n\times \Rm^n.$$
Because $h$ has positive definite Hessian,
the equation $\partial_{p_2}h(p_1,p_2)=0$ is non-singular and it defines
a smooth $m$-dimensional  manifold in $\Rm^n$  wich can also be described
parametrically  by the relation  $p_2=P_2(p_1)$
for some function $P_2:\Rm^m\lto \Rm^r$.
Therefore, the corresponding invariant manifold can be written in a parametric form 
as 
$$
\big\{(q_1,0,p_1,P_2(p_1));(q_1,p_1) \in \Tm^m\times \Rm^m\big\},
$$
it is a cylinder.
Moreover, this manifold is normally hyperbolic in the sense of \cite{HPS}.
It is necessary at this point to precise the terminology.
An open manifold will be called \textbf{weakly invariant} for a flow
if the vector field  is tangent at each point to this manifold. It will be called
\textbf{strongly invariant} if it contains the full orbit of each of its points.
A compact strongly invariant manifold is called normally hyperbolic if it 
is eventually  absolutely $1$-normally hyperbolic for the time-one flow in the sense of 
\cite{HPS}, definition 4.
\begin{defn}\label{def}
A weakly invariant open manifold $N$ (for some vector field $X$) is called normally hyperbolic
if there exists:
\begin{itemize}
\item A vector field $Y$ on a compact manifold
$M$.
\item An embedding  $i:U \lto M$ from  a neighborhood
$U$ of $N$ into $M$ which conjugates $X_{|U}$ and $Y_{|i(U)}$.
\item A normally hyperbolic strongly invariant compact manifold
$\tilde N$  in $M$
(for the vector field $Y$) such that $i(N)\subset \tilde N$.
\end{itemize} 
\end{defn}

Returning to the invariant cylinder of the averaged system,
we observe that the open sub-cylinder
$$
\big\{(q_1,0,p_1,P_2(p_1));(q_1,p_1) \in \Tm^m\times \Rm^m,
\|p_1\|<\delta \big\}, \quad \delta>0
$$
is a normally hyperbolic weakly (and even strongly) invariant
open sub-manifold for the averaged system in the sense of 
Definition \ref{def}.
From this observation, and from the fact that the full
system can be considered locally (near $p=p_0$) 
as a perturbation of the averaged system, one can prove the 
existence of a small normally hyperbolic weakly invariant
cylinder in the full system, this is  well understood.
This cylinders can also be seen 
as the center manifold of a "whiskered" 
(or partially hyperbolic) torus,
which is the continuation in the full systems of the 
invariant torus 
$$\big\{(t,q_1,0,p_0),\quad (t,q_1)\in \Tm\times \Tm^m
\big\}
$$
which exists in the
 averaged system. The name whiskered comes 
from the fact that this torus has hyperbolic normal directions,
this name  (as well as the corresponding object) was introduced by Arnold in \cite{A}.
 The existence of a whiskered torus in the original system was proved in 
\cite{Tre}, following earlier works on the persistence of partially hyperbolic KAM tori.
It is well understood, see for example  \cite{BT}
that such a torus must be contained in an invariant cylinder
which is normally hyperbolic. 
Proving the existence of whiskered tori involves KAM theory, which
is quite demanding in terms of regularity, while the existence of the invariant
cylinder relies on the softer theory of normal hyperbolicity.
The idea of embedding whiskered tori into a normally hyperbolic cylinder
and to use the  theory of normal hyperbolicity in the context of Arnold 
diffusion  is more recent
than the paper of Arnold. To the best of our knowledge, it appears first in 
Moeckel \cite{Mo}. 
It was then progressively understood  that normally hyperbolic
invariant cylinders can be used to produce diffusion even in the absence 
of whiskered tori. 

We described  two well-known methods allowing  to prove the existence
of small normally hyperbolic weakly invariant cylinders in
the full system for $\epsilon>0$. 
However, the size of the invariant cylinder that has been obtained 
in the literature is small, meaning that it converges to $0$
with $\epsilon$. Our point in the present paper is that a large 
normally hyperbolic weakly invariant cylinder actually  exists:

\begin{thm}\label{main}
Assume that $H$ is smooth (or at least $C^r$ for a sufficiently large $r$)
 and satisfies Hypothesis \ref{nondeg}.
Assume that $\omega$ is Diophantine, and fix $\kappa>0$.
Then there exists an  open ball $B\subset \Rm^m$ containing $p_1^0$,
a neighborhood $U$ of $0$ in $\Tm^r$,
 a positive number $\epsilon_0$
and,  for $\epsilon <\epsilon_0$ two $C^1$  functions
$$Q_2^{\epsilon}:\Tm\times \Tm^m\times B\lto U\subset  \Tm^r
\quad \text{ and }\quad 
P_2^{\epsilon}:\Tm\times \Tm^m\times B\lto \Rm^r
$$
such that the annulus
$$
A^{\e}=\big\{(t,q_1,Q_2^{\epsilon}(t,q_1,p_1),p_1,P_2^{\epsilon}(t,q_1,p_1))
, \quad (t,q_1,p_1)\in \Tm\times \Tm^m\times B\big\}
$$
is weakly invariant for (\ref{H})
(in the sense that the Hamiltonian vector field is tangent to it).
We have $P_2^{\epsilon}\lto P^0_2$ uniformly as $\epsilon\lto 0$,
where $P_2^0$ is the function $(t,q_1,p_1)\lmto P_2(p_1)$.
Moreover, we have $\|P_2^{\e}-P_2^0\|_{C^1}\leq \kappa$, and
 $\|Q_2^{\e}\|_{C^1}\leq \kappa /\epsilon$.
Each strongly invariant set of (\ref{H}) (in the sense that it contains the full orbit of each
of its points, for example, a whiskered torus)  contained in
the domain
$$
\mD^{\e}:=\Tm\times \Tm^m \times U\times B\times \{p_2\in \Rm^r : \|p_2\|\leq \e \}
$$
is contained in $A^{\e}$   for $\e < \e_0$.  The cylinder $A^{\e}$ is 
 normally 
hyperbolic  and  symplectic.
\end{thm}

The novelty  here is that the ball $B$
does not depend on $\epsilon$.
Easy examples show that we can't expect a fine  control 
of the asymptotic behavior of $Q^{\epsilon}_2$ in terms of the 
averaged system only except if we restrict to
smaller domains depending on $\epsilon$.
This asymptotic behavior  also depends on the 
averaged systems at other frequencies.
However, the very weak estimates we have are sufficient to 
describe  the restricted dynamics.
Let $A^{\e}_0\subset \Tm^n\times \Rm^n$ be the 
restriction of the invariant annulus to the section $\{t=0\}$, 
$$
A^{\e}_0=\big\{(q_1,Q_2^{\epsilon}(0,q_1,p_1),p_1,P_2^{\epsilon}(0,q_1,p_1))
, \quad (q_1,p_1)\in \Tm^m\times B\big\},
$$
and let $\phi:
A^{\e}_0\lto \Tm^n\times \Rm^n$
be the time-one flow of $H$ (which is
well-defined  on $A^{\e}_0$ when $\e$ is small enough).
Then $A^{\e}_0$ is somewhat invariant for  $\phi$
(although there are some difficulties near the boundary)
in a sense that will be given more precisely below. 
We  define the map $\Phi:\Tm^m\times B\lto \Tm^m\times \Rm^m$
 as the restriction of $\phi$ to $A^{\e}_0$
seen in coordinates $(q_1,p_1)$, more precisely
$$
\Phi(q_1,p_1)=(q_1,p_1)\circ \phi\big( q_1,Q_2^{\epsilon}(0,q_1,p_1),p_1,P_2^{\epsilon}(0,q_1,p_1)\big).
$$
Note that this map is well-defined on $ \Tm^m\times B$. 
Let us finally consider an open ball
 $B_0\subset \Rm^m$  which contains $p_0$
and whose closure is contained in $B$, and set
$$
A^{\e}_{00}=\big\{(q_1,Q_2^{\epsilon}(0,q_1,p_1),p_1,P_2^{\epsilon}(0,q_1,p_1))
, \quad (q_1,p_1)\in \Tm^m\times B_0\big\}.
$$

\begin{prop}\label{twist}
The map $\Phi$ is converging uniformly (when $\epsilon \lto 0$) on $\Tm^m\times B_0$ to the map
 $$
\Phi_0:\left( \begin{matrix}
q_1\\p_1
\end{matrix}\right)
\lmto 
\left( \begin{matrix}
q_1+\partial_{p_1}h\big(p_1,P_2(p_1)\big)
\\p_1
\end{matrix}
\right),
$$
which gives the unperturbed dynamics on the invariant cylinder of the averaged system.
Moreover, we have 
 $\phi(A^{\e}_{00})\subset A^{\e}_0$
when $\epsilon$ is small enough.
Finally, given $\eta>0$, we can choose the ball $B_0$ small enough so that the inequality
$$
\|d\Phi-d\Phi_0\|_{C^0}\leq \eta
$$
holds on $\Tm^m\times B_0$ when $\e$ is small enough.
\end{prop}

The  frequency map
$$
p_1\lmto \Omega_0(p_1):= \partial_{p_1}h\big(p_1,P_2(p_1)\big)
$$
has positive torsion  in the sense that 
$$
\partial_{p_1}\Omega_0=\partial^2_{p_1}h(p_1,P_2(p_1))
$$
is a positive definite symmetric matrix for all $p_1\in \Rm^m$.
As a consequence, when $\e$ is small enough, the restricted map $\Phi$ has positive torsion
in a neighborhood
(independent of $\e$) of $\Tm^m \times \{p_1^0\}$, in the sense that 
$$
\partial_{p_1}(q_1\circ \Phi)_{(q_1,p_1)}\rho_1\cdot \rho_1>0
\quad \forall \rho_1\in \Rm^m
$$
for all $q_1\in \Tm^m$ and $p_1\in B_0$, provided that $B_0$ has been chosen
small enough.
The map $\Phi$ is symplectic with respect to the symplectic
form obtained by restriction of the ambient symplectic form
to $A^{\e}_0$. 
It is part of the statement of Theorem \ref{main} that this form is 
non-degenerate on $A^{\e}_0$.
Note that this symplectic form is
not $dq_1\wedge dp_1$ in general.

In the case $m=1$,
(but for any dimension $n$) 
one can combine these results with existing techniques
on the
 \textit{a priori} unstable situation,
like the variational methods coming from Mather Theory 
(see \cite{Mather:93,fourier}),  developped for the  \textit{a priori} unstable situation in 
\cite{JAMS, CY, CY2} or more geometric methods like \cite{GR}
(The papers \cite{DLS:06,Tre} also treat the \textit{a priori} unstable situation, but it seems to me at first sight that they require too strong informations on the restricted dynamics to be applicable here).
One can then  hope to obtain, under additional non-degeneracy assumptions,
 the existence of restricted Arnold diffusion
in the following sense:
There exists $\delta>0$ and $\epsilon_0$ such that,
for each $\epsilon\in ]0,\epsilon_0[$
there exists an orbit $(q_{\e}(t), p_{\e}(t))$ with the following property:
The image $p_{\e}(\Rm)$ is not contained 
in any ball of radius $\delta$ in $\Rm^n$.
Once again, the key point here is that $\delta$
can be chosen independent of $\e$.
Specifying the needed "non-degeneracy assumptions"
 will require some further work, but I believe it will not require any method
beyond those which are  already available.

Of course, finding "global" Arnold diffusion,
as announced in \cite{M:A},
 that is orbits wondering in the 
whole phase space along different resonant lines (or far away along a
given resonant line)
 requires a specific study of relative resonances (when the restricted frequency
$\omega$ is resonant),
where the existence of normally hyperbolic invariant cylinders
can't be obtained by the method used in the present paper.

Let us close this introduction with a remark on uniqueness.
In general, there is no uniqueness statement for the normally invariant cylinder
we obtain. However, in the case $m=1$, we can obtain a stronger result:
Let $[p_1^-,p_1^+]\subset B\subset \Rm$ be an interval such that 
both $\Omega_0(p_1^-)$  and $\Omega_0(p_1^+)$ are Diophantine.
Then, there exists whiskered tori $\Tm^{\e}_-$ and $\Tm^{\e}_+$
 of dimension $2$ in $\Tm\times\Tm^n\times \Rm^n$
which are close to the unperturbed tori
$$T_-^0=\big\{\big(t,q_1,0,p_1^-,P_2(p_1^-)\big)
:\quad (t,q_1)\in \Tm\times \Tm\big\}
$$
and
$$T_+^0=\big\{\big(t,q_1,0,p_1^+,P_2(p_1^+\big)
:\quad (t,q_1)\in \Tm\times \Tm\big\}.
$$
The whiskered tori $\Tm^{\e}_{\pm}$ are  contained 
in the annulus $A^{\e}$. They bound a compact part $A^{\e}_{=}$ of $A^{\e}$
which is then stronly invariant in the  sense that it contains the full
orbit of each of its points. The annulus $A^{\e}_{=}$
is then unique in the sense that if $\tilde A^{\e}$ is another normally
hyperbolic cylinder given by Theorem \ref{main} (with the same domain $B$),
then it must contain $A^{\e}_{=}$. 
The cylinder $A^{\e}_{=}$ is a normally hyperbolic invariant cylinder
in the genuine sense.
 If the interval $[p^-,p^+]$ has been chosen small enough, then
the restricted map $\Phi:A^{\e}_=\lto A^{\e}_=$ is a $C^1$  area preserving twist map
(for the appropriate area form).
When $m>1$ one should not expect the same kind of properties,
 since Arnold diffusion may occur inside the invariant cylinder.

\section{Averaging}\label{average}
In order to apply averaging methods, it is 
easier to consider the extended phase space
$$(t,e,q,p)\in \Tm\times \Rm\times \Tm^n \times \Rm^n
$$
where the Hamiltonian flow can be seen
as the Hamiltonian flow of the autonomous Hamiltonian function
$$
\tilde H(t,e,q,p)=h(p)+e-\e^2 G(t,q,p)
$$
on one of its energy surfaces, for example $\tilde H=0$.
Then, we consider a smooth solution $f(t,q)$ 
of the Homological equation
$$
\partial_t f+ \partial_q f \cdot (\omega, 0) =G(t,q,p_0)-V(q_2).
$$
Such a solution exists because $\omega$ is Diophantine, as can be checked
easily by power series expansion. It is unique up to an additive constant.
We consider the smooth symplectic diffeomorphism
$$
\psi^{\e}: (t,e,q,p)\lmto (t,e+\e^2 \partial_tf(t,q),q, p+\e^2\partial_q f(t,q))
$$
and use the same notation for the diffeomorphism
$(t,q,p)\lmto (t,q, p+\e^2\partial_q f(t,q))$.
We have 
$$
\tilde H \circ \psi^{\e}
=
h(p)+e -\e^2 V(q_2) -\e^2R(t,q,p) +O(\e^4),
$$
where $R(t,q,p)=G(t,q,p)-G(t,q,p_0)$.
In other words, by the time-dependent symplectic change of coordinates 
$\psi^{\e}$, we have reduced the study of $H$ to the study of 
the time-dependent Hamiltonian
$$
H_1(t,q,p)=h(p)-\e^2 V(q_2)-\e^2 R(t,q,p) +O(\e^4)
$$
where $R=O(p-p_0)$.
As a consequence, Theorem \ref{main}
holds for $H$ if it holds for $H_1$.
More precisely,
assume that there exists an invariant cylinder
$$\mA^{\e} = 
(t,q_1,\mQ_2^{\epsilon}(t,q_1,p_1),p_1,\mP_2^{\epsilon}(t,q_1,p_1))
$$
for $H_1$, with $\|\mQ_2^{\epsilon}\|_{C^1}\leq \kappa/2\e$ and 
$\|\mP^{\e}_2-P^0_2\|_{C^1}\leq \kappa/2$.
Then the annulus $A^{\e}:=\psi^{\e}(\mA^{\e})$
is invariant for $H$. 
Since $\psi^{\e}$ is $\e^2$-close to the identity, while
$\|\mQ_2^{\epsilon}\|_{C^1}\leq \kappa/2\e$,
the annulus $A^{\e}$ has the form
$$
 A^{\e} = 
(t,q_1,Q_2^{\epsilon}(t,q_1,p_1),p_1,P_2^{\epsilon}(t,q_1,p_1))
$$
for $C^1$ functions $Q_2^{\e}$, $P^{\e}_2$ which satisfy
$\|Q_2^{\epsilon}\|_{C^1}\leq \kappa/\e$ and
$\|P^{\e}_2-P^0_2\|_{C^1}\leq \kappa$. 
We will prove that Theorem \ref{main} holds for $H_1$ in section 
\ref{proof}. We first expose some useful tools.

\section{Normally hyperbolic manifolds}\label{fix}

We shall now present a version of the classical theory of normally
hyperbolic manifolds adapted for our purpose.
On $\Rm^{n_z}\times \Rm^{n_x} \times \Rm^{n_y}$,
let us consider the time dependent  vector field
\begin{alignat*}{1}
\dot z&=Z(t,z,x,y)\\
\dot x &=A(z) x\\
\dot y & =-B(z) y.
\end{alignat*}
We assume that the  function
$$Z:\Rm\times \Rm^{n_z}\times \Rm^{n_x}\times \Rm^{n_y} \lto \Rm^{n_z}
$$
is $C^1$-bounded in the domain
\begin{equation}\tag{D}\label{D}
\Rm\times \Rm^{n_z}\times \{x \in \Rm^{n_x} : \|x\|< 1\}
\times \{y \in \Rm^{n_y} : \|y\|<1\},
\end{equation}
 and that the matrices $A$ and $B$ are $C^1$-bounded 
 functions of $z$.
 Moreover, we assume that there exists  constants
$a>b>0$ such that 
$$A(z)x\cdot x\geq a\|x\|^2 \quad ,\quad B(z)y\cdot y\geq a\|y\|^2
$$
for all $x,y,z$,
and such that 
$$
\|\partial_{(t,z)} Z(t,z,x,y)\|\leq b
$$
for all $(t,z,x,y)$ belonging to (\ref{D}).
We consider the perturbed vector field 
\begin{alignat*}{3}
\dot z&=Z(t,z,x,y)   &+&R_z(t,z,x,y)\\
\dot x &=A(z) x      &+&R_x(t,z,x,y)\\
\dot y & =-B(z) y   &+& R_x(t,z,x,y).
\end{alignat*}
where $R=(R_z,R_x,R_y)$ is seen as a small perturbation.

\begin{thm}\label{normal}
There exists $\e>0$ such that, when $\|R\|_{C^1}<\e$,
the maximal invariant set of the perturbed vector field
contained in the domain (\ref{D})
is a graph of the form 
$$
\big\{(t,z,X(t,z),Y(t,z)), \quad (t,z)\in \Rm\times \Rm^{n_z} \big\}
$$
where  $X$ and $Y$  are $C^1$ maps.
This graph is normally hyperbolic, and it is contained in 
the domain 
$$\Rm^{n_z}\times \{x \in \Rm^{n_x} : \|x\|\leq (2/a)\|R\|_{C^0}\}
\times \{y \in \Rm^{n_y} : \|y\|\leq (2/a)\|R\|_{C^0}\}.
$$
In other words, we have 
$$
\|(X,Y)\|_{C^0}\leq (2/a)\|R\|_{C^0}.
$$
The $C^1$ norm of $(X,Y)$ is converging to zero
when the $C^1$ norm of the perturbation converges to zero.
\end{thm}

\proof
The invariant space $\Rm^{n_z}$ is normally hyperbolic in the sense of
\cite{F:71, HPS}.
As a consequence, the standard theory applies and implies the existence of functions 
$X$ and $Y$ such that the graph
 $
(t,z,X(t,z),Y(t,z))
$
is invariant, normally hyperbolic, and contained in (\ref{D}).
Note that we are slightly outside of the hypotheses of the statements in 
\cite{HPS} because our unperturbed manifold is not compact.
However, the results  actually  depend on uniform estimates
rather than on compactness (see \cite{DLS}, Appendix B, for example, see also \cite{C}), 
and we assumed such uniform estimates.

Let us now prove the estimate on $(X,Y)$.
We have the inequality 
$$\dot x \cdot x\geq a\|x\|^2 +x\cdot R_x\geq a\|x\|(\|x\|-\|R_x\|_{C^0}/a)
$$
which implies that 
$$
\dot x\cdot x \geq \|x\|\|R_x\|_{C^0}
$$
if 
$$2\|R_x\|_{C^0}/a\leq \|x\|\leq 1,
$$
hence this domain can't intersect the invariant graph.
Similar considerations show that the domain
$2\|R_y\|_{C^0}/a\leq \|y\|\leq 1$
can't intersect the graph.
\qed

\section{Hyperbolic Linear System}

Let us consider the linear Hamiltonian system on $\Rm^n\times \Rm^n$
generated by the Hamiltonian
$$
H(q,p)=\frac{1}{2}\langle Bp,p\rangle-\frac{1}{2}\langle Aq,q\rangle,
$$
where both $A$ and $B$ are positive definite symmetric matrices.
We recall that this system can be reduced to
$$
G(x,y)=\langle Dx,y\rangle,
$$
where $D$ is a positive definite symmetric matrix,
by a linear symplectic change of variables $(q,p)\lto (x,y)$.
In order to do so, we consider the symmetric positive definite matrix
$$
L:=\big(A^{-1/2}(A^{1/2}BA^{1/2})^{1/2}A^{-1/2}\big)^{1/2},
$$
which is the only symmetric and positive definite solution of the equation 
$L^2AL^2=B$.
Considering the change of variables
$$ x=\frac{1}{\sqrt{2}}(Lp+L^{-1}q)\quad;\quad
y=\frac{1}{\sqrt{2}}(Lp-L^{-1}q)$$
or equivalently
$$
q=\frac{1}{\sqrt{2}}L(x-y) \quad;
\quad
p=\frac{1}{\sqrt{2}}L^{-1}(x+y),
$$
an elementary calculation shows that we obtain the desired form
for the Hamiltonian in coordinates $(x,y)$, with 
$$
D=LAL=L^{-1}BL^{-1}.
$$
As a consequence, the equations of motions in the new variables 
take the block-diagonal form
$$\dot x=Dx\quad ; \quad \dot y=-Dy.
$$
In the original coordinates $(q,p)$ the stable space (which is the space $x=0$)
is the space $\{(q,-L^2q), q\in \Rm^n\}$ while the unstable space is 
$\{(q,L^2q), q\in \Rm^n\}$.

\section{Proof of Theorem \ref{main}}\label{proof}
We now prove Theorem \ref{main} for the Hamiltonian
$$
H_1(t,q,p)=h(p)-\epsilon^2 V(q_2) -\epsilon^2R(t,q,p) +O(\epsilon^{2+\gamma}),
$$
where $R=O(p-p_0)$ and $\gamma>0$ ($\gamma=2$ in our situation). We assume that Hypothesis \ref{nondeg} holds.
 We lift all the angular variables to the 
universal covering, and see $H_1$ as a Hamiltonian of the variables 
$$
(t,q,p)=(t,q_1,q_2,p_1,p_2)\in 
\Rm\times \Rm^m\times\Rm^r \times \Rm^m \times \Rm^r
$$
which is one-periodic in $t,q$. We assume that $p_0=0$.

We will need  some notations.
We set $A:= \partial^2V(0)$, it is a symmetric positive definite matrix.
We will denote by $B(p_1)$ a matrix which depends smoothly on $p_1$, 
is uniformly positive definite, is constant outside of a neighborhood
of $p_1=0$ in $\Rm^m$, and coincides with $\partial^2_{p_2}h(p_1,P_2(p_1))$
in a neighborhood of $p_1=0$. We will denote by $\tilde P_2(p_1)$
a compactly supported smooth function $\tilde P_2:\Rm^m\lto \Rm^r$ which coincides 
with $P_2$ around $p_1=0$. 
Finally, we will denote by $h_0(p_1)$ a smooth compactly supported function
which is equal to $h(p_1,P_2(p_1))$ around $p_1=0$.

It is useful to introduce two new positive parameters $\alpha$ and $\delta$.
We always assume that 
$$
0<\e <\delta<\alpha <1.
$$
In the sequel, we shall chose $\alpha$ small, then $\delta$ small with respect
to $\alpha$, and work with $\e$ small enough with respect to $\alpha$ and $\delta$. The parameter $\delta$ represents the size of
the normally hyperbolic cylinder we intend to find. 
We will denote by $\underline \chi$ a smooth function of its arguments 
which may depend (in an unexplicited way) on the parameters
$\e, \delta$, but which is $C^2$-bounded, uniformly in 
$\e, \delta$. The notation $\chi$ will be used in a similar way
when only $C^1$ bounds are assumed.

\begin{lem}
There exists a smooth Hamiltonian function $H_2(t,q,p)$
(which depends on the parameters $\e, \delta$)
of the form
\begin{align*}
H_2&=h_0(p_1)+\2B(p_1)\cdot (p_2-\tilde P_2(p_1))^2
-\frac{\e^2}{2}A\cdot q_2^2 \\
	&+\e^3\underline\chi \big(p_1,(p_2-\tilde P_2(p_1))/\e\big)
	+\e^2\delta^{3/2}\underline\chi(q_2/\sqrt \delta)
	+\e^2\delta\underline \chi (t,q,p/\delta)
	+\e^{2+\gamma}\underline\chi(t,q,p)
\end{align*}
which coincides with $H_1$
on the domain
$$\big\{\|q_2\|\leq \sqrt \delta, \|p_1\|\leq \delta, \|p_2-P_2(p_1)\|\leq \e
\big\}.
$$
\end{lem}
\proof
Let us expand the function $h$ with respect to $p_2$ at the point $P_2(p_1)$:
$$
h(p_1,p_2)
=h(p_1,P_2(p_1))
+\2  \partial_{p_2}^2h(p_1,P_2(p_1))\cdot (p_2-P_2(p_1))^2
+\underline S(p)\cdot (p_1-P_2(p_1))^3
$$
 where $\underline  S(p) $
is a $3$-linear form on $\Rm^r$ depending smoothly on $p$.
We consider a $3$-form $S(p)$ which depends smoothly on $p$, is compactly
supported, and is equal to $\underline S(p)$ near $p=0$.
Let $ i:\Rm^k\lto \Rm^k$ (for any $k$) be a compactly
supported smooth map which is equal to the identity on the unit ball.
Then the function
\begin{align*}
&h_0(p_1)+\2 B(p_1)\cdot (p_1-\tilde P_2(p_1))^2+\e^3S(p)\cdot
\big(i\big[(p_2-\tilde P_2(p_1))/\e\big]
\big)^3\\
=&h_0(p_1)+\2 B(p_1)\cdot (p_1-\tilde P_2(p_1))^2+
\e^3\underline \chi\big(p_1,(p_2-\tilde P_2(p_1))/\e\big)
\end{align*}
is equal to $h$ if $p$ belongs to a given neighborhood of $0$
(independant of $\e, \delta$) and satisfies
$\|p_2-\tilde P_2(p_1)\|\leq \e$.
Similarly,
we write
$$
V(q_2)=\2 A\cdot q_2^2 +W(q_2)\cdot q_2^3
$$
for some $3$-linear form $W(q_2)$.
It is equal to
$$
\2 A\cdot q_2^2 +\delta^{3/2}W(q_2)\cdot (i(q_2/\sqrt\delta))^3=
\2 A\cdot q_2^2 +\delta^{3/2}\underline \chi \big(q_2/\sqrt\delta\big)
$$
on $\{\|q_2\|\leq \sqrt \delta\}$.
Finally, we observe that the function 
$R(t,q,p)$ can be written in the form
$$
R(t,q,p)=L(t,q,p)\cdot p
$$
and is equal to the function
$$
\delta L(t,q,p)\cdot i(p/\delta)=\delta \underline \chi(t,q,p/\delta)
$$
on $\{\|p\|\leq \delta\}$.
Collecting all terms proves the Lemma.
\qed

We will now prove the existence of a normally hyperbolic invariant graph
for $H_2$ contained in the region 
$$\{ \|q_2\|\leq \sqrt\delta, 
\|p_2-\tilde P_2(p_1)\|\leq \e  \}
$$
Its intersection with $\{\|p_1\|< \delta\}$ will give a weakly invariant
manifold for $H_1$
(meaning that the Hamiltonian vector field of $H_1$ is tangent to it).
In order to simplify the following equations,
we set 
$$
h_2(p):= h_0(p_1)+\2 B(p_1)\cdot (p_2-\tilde P_2(p_1))^2.
$$
The Hamiltonian vector field of $H_2$ can be written
\begin{alignat*}{3}
\dot q_1& = \partial_{p_1}h_2(p)
	&+&\e^2 \chi \big(p_1,(p_2-\tilde P_2(p_1))/\e\big)+\e^2 \chi(t,q,p/\delta) \\
\dot p_1& =0 &+& \epsilon^2 \delta\chi(t,q,p/\delta)\\
\dot q_2& = B(p_1)(p_2-\tilde P_2(p_1))&+&
	\epsilon^2  \chi(p_1,(p_2-\tilde  P_2(p_1))/\epsilon)+\e^{2}\chi(t,q,p)\\
\dot p_2& =\e^2 Aq_2 &+&\e ^2\delta  \chi(q_2/\sqrt \delta)+
 	\e^2 \delta \chi(t,q_1,q_2,p/\delta)
\end{alignat*}
recalling the convention that $\chi(.)$ always denotes
a $C^1$ function of its arguments, depending on $\e$ and $\delta$,
but bounded in $C^1$ independently of $\delta $ and $\epsilon$.
Motivated by section \ref{fix}, we set 
$$
L(p_1)=\big(A^{-1/2}(A^{1/2}B(p_1)A^{1/2})^{1/2}A^{-1/2}\big)^{1/2},
$$
and  perform the change  of variables
$(t,q_1,p_1,q_2,p_2)\lto (\tau,\theta, r, x,y)$ given by:
\begin{align*}
\tau =\e t,\quad \theta = \e\alpha q_1,&\quad r=p_1,\\
 x=L(p_1)(p_2-\tilde P_2(p_1))+\e L^{-1}(p_1)q_2,&\quad
y=L(p_1)(p_2-\tilde P_2(p_1))-\e L^{-1}(p_1)q_2,
\end{align*}
recalling that $\alpha$ is a fixed positive parameter.
Equivalently, this can be written
$$
t=\tau/\e ,\quad q_1=\theta/\e\alpha,\quad p_1=r,\quad
 q_2=L(r)(x-y)/2\e,\quad p_2=\tilde P_2(r)+L^{-1}(r)(x+y)/2.
$$
In the new coordinates, the principal part of the vector field
takes the form (denoting $\acute f$ for $df/d\tau$) 
$$
\acute \theta =\alpha \Omega(r,x,y),\quad
\acute r=0, \quad
\acute x= D(r)x,\quad
\acute y=- D(r)y,
$$
with 
$$\Omega(r,x,y):=\partial_{p_1}h_2
\big(r,\tilde P_2(r)+L^{-1}(r)(x+y)/2\big)
$$
and
$$
D(r):=L(r)AL(r)=L^{-1}(r)B(r)L^{-1}(r).
$$
The equality above holds because $L(r)$ solves the equation 
$L^2(r)AL^2(r)=B(r)$.
Let us detail the calculations leading to the expressions
of $\acute x:= dx/d\tau$ (the calculation for $\acute y$ is similar):
\begin{align*}
\e \acute x=\dot x=&
L(p_1)\big( \dot p_2-\partial_{p_1}\tilde P_2\cdot \dot p_1\big)
+\e L^{-1}(p_1)\dot q_2
+ \big(\partial_{p_1} L\cdot \dot p_1 \big)
\big(p_2-\tilde P_2(p_1)\big)
+ 
\e\big( \partial_{p_1}(L^{-1})\cdot \dot p_1\big)
 q_2\\
=& \e^2 L(p_1) A q_2+ \e L^{-1}(p_1)B(p_1)(p_2-\tilde P_2(p_1))\\
+& \e^2\delta \chi(t,q,p/\delta,x,y)+\e^3\chi(p_1,(x+y)/\e)+
\e^2\delta\chi(q_2/\sqrt{\delta})+\e^{2+\gamma}\chi(t,q,p,x,y)
\\
=& \e L(r)A L(r)(x-y)/2+\e L^{-1}(r)B(r)L^{-1}(r)(x+y)/2\\
+& \e^2\delta \chi(\tau/\e,\theta/\e,r/\delta ,x/\delta, y/\delta,x/\e,y/\e)+
\e^2\delta \chi(r/\sqrt{\delta},x/\sqrt{\delta}\e, y/\sqrt{\delta}\e)
\\
=& \e D(r) x 
+\e^2\delta \chi(\tau/\e,\theta/\e,r/\delta ,x/\delta, y/\delta,x/\e,y/\e)+
\e^2\delta \chi(r/\sqrt{\delta},x/\sqrt{\delta}\e, y/\sqrt{\delta}\e).
\end{align*}
The function $\Omega (r,x,y)$ is $C^1$-bounded on 
$$
\big\{(r,x,y), \quad\|x\|\leq 1, \|y\|\leq 1\big\}.
$$
We can choose $\alpha<1$ once and for all in order that the principal part 
of the vector field satisfies the hypotheses of Theorem \ref{normal}.
The full vector field can be written in the new coordinates,
(with the notation  $\acute f:= df/d\tau$):
\begin{alignat*}{3}
\acute \theta & =\alpha  \Omega(r,x,y)
	&+&\e^2 \chi(\tau/\e,\theta/\alpha\e,r,x/\e,y/\e) \\
\acute r& =0 &+&
 	\e \delta\chi(\tau/\e, \theta/\alpha\e, r/\delta,x/\delta, y/\delta, x/\e, y/\e)\\
\acute x& =D(r)x&+&
	\e\delta \chi(\tau/\e,\theta/\alpha\e,r/\delta ,x/\delta, y/\delta,x/\e,y/\e)+
	\e\delta \chi(r/\sqrt{\delta},x/\sqrt{\delta}\e, y/\sqrt{\delta}\e)\\
\acute y& =-D(r)y&+&
	\e\delta \chi(\tau/\e,\theta/\alpha\e,r/\delta ,x/\delta, y/\delta,x/\e,y/\e)+
	\e\delta \chi(r/\sqrt{\delta},x/\sqrt{\delta}\e, y/\sqrt{\delta}\e).
\end{alignat*}
In this expression, 
we observe that the uniform norm of the perturbation is 
$O(\e\delta)$ while the $C^1$ norm is 
$O(\sqrt \delta)$  (recall that $0<\epsilon< \delta < 1$).
We can apply Theorem \ref{normal}
and find a unique bounded normally hyperbolic  invariant graph
$$
(\tau,\theta, X(\tau,\theta,r),r,Y(\tau, \theta,r)).
$$ 
Moreover Theorem \ref{normal} also implies that 
$$\|(X,Y)\|_{C^0} \leq C\e \delta.$$
Since the invariant graph we have obtained is the maximal invariant set
contained in the domain $\{\|x\|\leq 1, \|y\|\leq 1\}$, and since the vector field
is $\e$-periodic in $t$ and $\alpha\e$-periodic in $q_1$, we conclude that the functions 
$X$ and $Y$ are $\e$-periodic in $t$ and $\alpha\e$-periodic in $q_1$.
In the initial coordinates, we have an invariant graph
$$
(t, q_1, Q_2^{\e}(t,q_1,p_1), p_1, P_2^{\e}(t,q_1,p_1))
$$
with
$$
Q^{\e}_2(t,q_1,p_1)= L(p_1) \big( 
X(\epsilon t, \e q_1, p_1)-Y(\e t, \e q_1, p_1)\big)/2\e 
$$
and 
$$
P^{\e}_2(t,q_1,p_1)=\tilde P_2(p_1)+
 L^{-1}(p_1) \big( 
X(\epsilon t, \e q_1, p_1)+Y(\e t, \e q_1, p_1)\big)/2.
$$
The functions $Q^{\e}_2$ and $P^{\e}_2$ are $1$-periodic in $(t,q_1)$.
The invariant graph we have obtained is normally hyperbolic 
for the flow of $H_2$, and its strong stable and strong unstable
directions  have 
the same dimension $r$. It follows from general results on partial hyperbolicity
in a symplectic context (see \textit{e. g.} \cite{LMS}, Proposition 1.8.3
\footnote{ In this text, the equality of the dimensions
of
the stable and unstable directions  (that obviously holds here) 
is stated  as a conclusion, although it should be taken as an assumption.})
that it is a symplectic manifold. This means that the restriction 
to the invariant graph of the ambient symplectic form
is a symplectic form.
Observing that 
$$\|Q_2^{\e}\|_{C^0}\leq C\delta, \quad
\|P_2^{\e}\|_{C^0}\leq C\e \delta,
$$
we infer  that the annulus
$$
\big\{
(t, q_1, Q_2^{\e}(t,q_1,p_1), p_1, P_2^{\e}(t,q_1,p_1)): 
t\in \Tm, q_1\in \Tm^m, p_1\in \Rm^m, \|p_1\|< \delta 
\big\}
\subset \Tm \times \Tm^n \times \Rm^n
$$
is contained in the domain 
$$
\{\|q_2\|\leq \sqrt \delta, \|p_1\|\leq \delta, \|p_2-\tilde P_2(p_1)\|\leq \e\}
$$
where $H_2=H_1$, provided $\delta$ has been chosen small enough. It is thus a weakly  invariant cylinder for $H_1$
\textit{i.e.} the extended Hamiltonian vector field of $H_1$
on $\Tm\times \Tm^n\times \Rm^n$ is tangent to this annulus at each point.
Orbits may still exit  from the cylinder through its boundary.
We finish with the estimates on the $C^1$ norms.
Since the $C^1$ size of the perturbation is $O(\sqrt{\delta})$,
we can make it as small as we want by chosing $\delta$
small.
We can thus assume that 
$\|(X,Y)\|_{C^1}
$
is small, and this implies the desired $C^1$ estimates on 
$P^{\e}_2$ and $Q^{\e}_2$.
We have proved Theorem \ref{main} for $H_1$, we conclude from Section \ref{average}
that Theorem \ref{main} holds for $H$.

\section{Proof of Proposition \ref{twist}}
Let $(q_1,p_1)$ be given in $\Tm^m\times B$,
and let $(q_1(t), q_2(t),p_1(t),p_2(t))$
be the orbit (under $H$)  of the point 
$$
\big(q_1, Q_2^{\e}(0,q_1,p_1),p_1, P_2^{\e}(0,q_2,p_2)\big).
$$
We have the Hamilton  equations
\begin{align*}
\dot q_1(t)&=\partial_{p_1}H(t,q_1(t),q_2(t),p_1(t),p_2(t))\\
\dot p_1(t)&=-\partial_{q_1}H(t,q_1(t),q_2(t),p_1(t),p_2(t)).
\end {align*}
They imply that $\dot p_1=O(\e^2)$, 
and we conclude that $p_1(t)\in B$ for all $t\in[0,1]$ if $p_1\in B_0$,
provided $\e$ is small enough.
This implies the inclusion 
$$
\phi(A^{\e}_{00})\subset A^{\e}_0,
$$
and it also implies that 
$$
\big(q_1(t), q_2(t),p_1(t),p_2(t)\big)=
\big(q_1(t), Q_2^{\e}(t,q_1(t),p_1(t)),
p_1(t), P_2^{\e}(t,q_1(t),p_1(t))\big)
$$
for each $t\in [0,1]$.
The Hamilton equations then take the form
\begin{alignat*}{2}
\dot q_1(t)&=\partial_{p_1}h\big(p_1(t),P_2^{\e}(t,q_1(t),p_1(t))\big)
		&-\e ^2 \partial_{p_1}G\big(t,q_1(t),Q_2^{\e}(t,q_1(t),p_1(t)),
		P_2^{\e}(t,q_1(t),p_1(t))\big)
\\
\dot p_1(t)&=&+\epsilon^2\partial_{q_1}G\big(t,q_1(t),Q_2^{\e}(t,q_1(t),p_1(t)),
		P_2^{\e}(t,q_1(t),p_1(t))\big).
\end{alignat*}
The map $\Phi$ is thus the time-one flow
of the vector field 
$$
\left(\begin{matrix}q_1\\
p_1\end{matrix}
\right)
\lmto 
\left(\begin{matrix}\partial_{p_1}h\big(p_1,P_2^{\e}(t,q_1,p_1)\big)-\e^2
 \partial_{p_1}G\big(t,q_1,Q_2^{\e}(t,q_1,p_1),
		P_2^{\e}(t,q_1,p_1)\big)
\\
\epsilon^2\partial_{q_1}G\big(t,q_1,Q_2^{\e}(t,q_1,p_1),
		P_2^{\e}(t,q_1,p_1)\big)
\end{matrix}
\right)
$$
which converges uniformly 
to the vector field
$$
\left(\begin{matrix}q_1\\
p_1\end{matrix}
\right)
\lmto 
\left(\begin{matrix}\partial_{p_1}h\big(p_1,P_2(p_1)\big)
\\
0
\end{matrix}
\right)
$$
when $\e \lto 0$ on $\Tm^m\times B$.
We conclude that $\Phi$ 
is converging uniformly to $\Phi_0$ (as defined in Proposition \ref{twist}).
Moreover, we see that the $C^1$ distance between these two vector fields
is $O(\kappa)$ ($\kappa$ is a parameter introduced in the
 statement of Theorem  \ref{main}),
so it can be made arbitrarily small by taking $B_0$ small enough.
The same statement then holds for the time-one flows $\Phi$ and $\Phi_0$.
\begin{small}

\end{small}

\begin{thebibliography}{99}
\bibitem{A}
\textsc{V. I. Arnold},
 Instability of dynamical systems with several degrees of freedom. 
\textit{Sov. Math. Doklady}, 
\textbf{5} (1964), 581--585.


\bibitem{fourier}
\textsc{P. Bernard}:  Connecting orbits of time dependent 
Lagrangian systems.
 \textit{Ann. Inst. Fourier} {\bf  52} (2002), 1533--1568.

\bibitem{JAMS}
\textsc{P. Bernard}: The dynamics of pseudographs
in convex Hamiltonian systems, Journ. AMS,
\textbf{21} No. 3 (2008) 625--669.

\bibitem{BT} 
\textsc{S.V. Bolotin, D.V. Treschev}:
{Remarks on the definition of hyperbolic tori of Hamiltonian systems}
 \textit{Regular and Chaotic dynamics},
  \textbf{5} (2000),
 no. {4},
  {401--412}.
 

\bibitem{C}
\textsc{M. Chaperon},
Stable manifolds and the Perron-Irwin method,
\textit{Erg. Th. Dyn. Sys.}, \textbf{24}  (2004), 1359-1394.  

\bibitem{CY}
\textsc{C.-Q. Cheng, J. Yan},
 Existence of diffusion orbits in a priori unstable Hamiltonian systems, 
 \textit{ J. Differential Geom.}
 \textbf{67} (2004), no. 3, 457--517.


\bibitem{CY2}
\textsc{C.-Q. Cheng, J. Yan},
 Arnold Diffusion in Hamiltonian systems:  the a priori unstable case,
 \textit{ J. Differential Geom.}
\textbf{82}, (2009),  no. 2, 229-277.


\bibitem{F:71}
\textsc {N. Fenichel}:{Persistence and smoothness of invariant manifolds for flows},
\textit{ Indiana Univ. Math. J.},
  \textbf{21}, (1971), 
{193--226}.

\bibitem{HPS}
 \textsc{M.W. Hirsch, C.C. Pugh, M. Shub}
\textit{Invariant manifolds},
Lecture notes in Math.  
 {Springer Berlin, New York}, (1977).


\bibitem{DLS:06}
\textsc{A. Delshams, R. de la Llave, T. M. Seara :}
\textit{A Geometric Mechanism for diffusion in Hamiltonian Systems Overcoming the 
Large Gap Problem: Heuristics  and Rigorous Verification on a Model,}
Mem. A.M.S. \textbf{179} (2006), no 844.



\bibitem{DLS}
\textsc{A. Delshams, R. de la Llave, T. M. Seara :}
Orbits of unbounded energy in quasi-periodic perturbations of geodesic flows,
\textit{Adv. in Math.} \textbf{202} (2006) 64-188.


\bibitem{LMS}
\textsc{P. Lochack, J. P. Marco, D. Sauzin :}
\textit{On the Splitting of Invariant Manifolds in 
Multidimensional Near-Integrable Hamiltonian Systems},
Mem. A.M.S. \textbf{163} (2003) no. 775.



\bibitem{GR}
\textsc{M. Gidea, C. Robinson}:
Obstruction argument for transition chains of 
Tori interspersed with gaps, preprint.


\bibitem{Mather:93}
\textsc{J. N. Mather}:
Variational construction of connecting orbits, 
Ann. Inst. Fourier, \textbf{43} (1993), 1349-1368.


\bibitem{M:A}
\textsc{J. N. Mather}: Arnold diffusion: announcement of results,
  J. Math. Sci. (N. Y.)  \textbf{124}  (2004),  no. 5, 5275--5289.


\bibitem{Mo}
\textsc{R. Moeckel}: Transition Tori in the Five-Body Problem,  JDE \textbf{129}
(1996), 290-314.


\bibitem{T}
\textsc{ D. Treschev}:
Hyperbolic tori and asymptotic surfaces in Hamiltonian systems
 Russ. J. Math. Phys, \textbf{2} (1994) no. 1, 93-110.


 \bibitem{Tre}
\textsc{D. Treschev},
 Evolution of slow variables in a priori unstable Hamiltonian systems. \textit{Nonlinearity}
 \textbf{ 17}  (2004), no. 5, 1803--1841. 


\end{thebibliography}
\end{document}